\documentstyle{amsppt}
\input bull-ppt
\keyedby{rosenthal/lic} 

\let\sl\it

\define\varep{\varepsilon}
\define\complex{{\Bbb C}}
\define\real{{\Bbb R}}
\define\dfeq{\buildrel {\text{\rm df}}\over =} 
\define\un#1{{\underline{#1}}}
%\define\olim{\mathop{\overline{\text{\rm lim}}}}
\predefine\preolim{\olim}
\redefine\olim{\mathop{\overline{\text{\rm lim}}}}
\define\osc{\operatorname{osc}}
%\define\Re{\text{\rm Re}}
\predefine\preRe{\Re}
\redefine\Re{\text{\rm Re}}
\define\Se{\text{\rm Se}}
\define\cc{\text{\rm (c.c.)}}
%\define\ss{\text{\rm (s.s.)}}
\predefine\press{\ss}
\redefine\ss{\text{\rm (s.s.)}}
\define\s{\text{\rm (s)}}
%\define\c{\text{\rm (c)}}
\predefine\underc{\c}
\redefine\c{\text{\rm (c)}}
\define\U{{\Cal U}} 
\define\V{{\Cal V}}

\topmatter 
\cvol{30}
\cvolyear{1994}
\cmonth{April}
\cyear{1994}
\cvolno{2}
\cpgs{227-233}
\ratitle
\title A subsequence principle characterizing\\
Banach spaces containing $c_0$\endtitle 
\shorttitle{A Subsequence Principle Characterizing Spaces 
Containing $c_0$}
\author Haskell Rosenthal\endauthor 
\thanks This research was 
partially supported by NSF DMS-8903197 and TARP 
235\endthanks%.\endthanks
\thanks The results stated here were presented on January 
16, 1993,
in the Special Session on Banach Space Theory, AMS meeting 
no. 878,
San Antonio, Texas\endthanks
\address Department of Mathematics,
The University of Texas at Austin,
Austin, Texas 78712-1082\endaddress%\endaddress
\ml%\email 
rosenthl\@math.utexas.edu\endml%email
%\endsubjclass
Weakly sequentially complete dual, convex block basis, the 
$\ell^1$-theorem, differences of semi-continuous 
functions\endkeywords 
\date July 13, 1993\enddate
\abstract 
The notion of a strongly summing sequence is introduced. 
Such a sequence is weak-Cauchy, a basis for its closed 
linear span, and has 
the crucial property that the dual of this span is not 
weakly sequentially 
complete. The main result is: 
%\vskip1pt
%{\narrower
\proclaim{Theorem}  Every non-trivial weak-Cauchy sequence 
in a \RM(real or 
complex\/\RM) Banach space has either a strongly summing 
sequence or a 
convex block basis equivalent to the summing basis.
\endproclaim%\smallskip}
%\vskip1pt 
%\noindent 
(A weak-Cauchy sequence is called {\it non-trivial\/} if 
it is 
{\it non-weakly convergent\/}.) The following 
characterization of spaces 
containing $c_0$ is thus obtained, in the spirit of the 
author's 1974 
subsequence principle characterizing Banach spaces 
containing $\ell^1$. 
%\vskip1pt
%{\narrower
\proclaim{Corollary 1}  A Banach space $B$ contains no 
isomorph  of $c_0$ 
if and only if every non-trivial weak-Cauchy sequence in 
$B$ has a strongly 
summing subsequence.\endproclaim%\smallskip}
%\vskip1pt
%\noindent 
Combining the $c_0$-and $\ell^1$-theorems, one obtains
%\vskip1pt 
%{\narrower
\proclaim{Corollary 2}  If $B$ is a non-reflexive Banach 
space such that 
$X^*$ is weakly sequentially complete for all linear 
subspaces $X$ of $B$, 
then $c_0$ embeds in $B$.\endproclaim%\smallskip}
\endabstract 
\endtopmatter 

%\rightheadtext{A subsequence principle characterizing 
%spaces containing $c_0$}
\document

\heading%\head 
1. Introduction\endheading%.\endhead

When does a general Banach space contain one of the 
classical sequence 
spaces? In particular, when does it contain one of the 
non-reflexive ones, 
$\ell^1$ (the space of absolutely summable sequences) or 
$c_0$ (the space 
of sequences vanishing at infinity)? 

In 1974, the following subsequence dichotomy was 
established by the author 
for real scalars \cite{R1} and refined by L.~E.~Dor to 
cover the case of 
complex scalars \cite{Do} (cf.\ also \cite{R2} for a 
general exposition). 

\proclaim{Theorem 1.0} 
Every bounded sequence in a real or complex Banach space 
has either 
a weak-Cauchy subsequence or a subsequence equivalent to 
the standard 
$\ell^1$-basis. 
\endproclaim 

The following characterization is an immediate consequence. 

\proclaim{Corollary 1 \cite{R1}}  
A real or complex Banach space $B$ contains no isomorph of 
$\ell^1$ if and 
only if every bounded sequence in $B$ has a weak-Cauchy 
subsequence. 
\endproclaim 

Using standard functional analysis, one easily obtains 

\proclaim{Corollary 2} 
If $B$ is a non-reflexive Banach space so that $B$ is 
weakly 
sequentially complete, then $\ell^1$ embeds in $B$. 
\endproclaim 

To obtain analogous results to characterize spaces 
containing $c_0$, 
 the following concept is introduced: 

\dfn{Definition}  
A sequence $(b_j)$ in a Banach space is called {\it 
strongly summing\/}, or 
{\rm(}s.s.{\rm)}, if $(b_j)$ is a weak-Cauchy basic 
sequence so that whenever scalars 
$(c_j)$ satisfy $\sup_n \|\sum_{j=1}^n c_j b_j\|<\infty$,  
the 
series $\sum c_j$ converges. 
\enddfn

A simple permanence property: let $(b_j)$ be an (s.s.) 
basis for a Banach 
space $B$, and let $(b_j^*)$ be its biorthogonal 
functionals in $B^*$, i.e., 
$b_j^* (b_i) = \delta_{ij}$ for all $i$ and $j$.  
Then $(\sum_{j=1}^n b_j^*)_{n=1}^\infty$ is a non-trivial 
weak-Cauchy 
sequence in $B^*$; hence $B^*$ {\it fails\/} to be weakly 
sequentially 
complete. (A weak-Cauchy sequence is called {\it 
non-trivial\/} if it is 
{\it non-weakly convergent\/}.) 

The following is the main result of this announcement. 

\proclaim{Theorem 1.1} 
Every non-trivial weak-Cauchy sequence in a {\rm(}real or 
complex\/{\rm)} 
Banach space 
has either a strongly summing subsequence or a convex 
block basis 
equivalent to the summing basis. 
\endproclaim 

\demo{Remark} 
The two alternatives of Theorem 1.1 are easily seen to be 
mutually exclusive. 
\enddemo 

Recall that the {\it summing basis\/} denotes the unit 
vectors basis for 
the space of all converging series of scalars, denoted 
$\Se$, endowed with 
the norm $\|(c_j)\|_{\Se} = \sup_n |\sum_{i=1}^n c_i|$. 
It is a standard, simple result that $\Se$ is isomorphic 
to $c_0$, 
and hence we obtain the following analogue of Corollary~1. 

\proclaim{Corollary 1.1}  
A Banach space $B$ contains no isomorph of $c_0$ if and 
only if every 
non-trivial weak-Cauchy sequence in $B$ has an $\ss$ 
subsequence. 
\endproclaim 

Combining the $c_0$- and $\ell^1$-theorems, we obtain the 
following result, 
analogous to Corollary~2. 

\proclaim{Corollary 1.2} 
If $B$ is a non-reflexive Banach space such that $X^*$ is 
weakly 
sequentially complete for all linear subspaces $X$ of $B$, 
then $c_0$ 
embeds in $B$. 
\endproclaim 

We note that the hereditary hypothesis in Corollary~1.2 is 
crucial. 
Indeed, J.~Bourgain and F.~Delbaen construct in 
\cite{Bo-De} a Banach 
space $B$ containing  {\it no\/} isomorph of $c_0$, such 
that $B^*$ is 
isomorphic to $\ell^1$ (so of course $B^*$ is weakly 
sequentially complete). 

To prove the main result, Theorem 1.1, we develop various 
permanence 
properties of strongly summing sequences and 
then employ some transfinite invariants for general 
discontinuous 
functions introduced by A.~S.~Kechris and A. Louveau 
\cite{AL}. 
These are used to characterize (complex) differences of 
bounded 
semi-continuous functions, which are essentially involved 
in the proof. 
The core of the argument is then a real-variables theorem 
concerning a 
subsequence refinement principle for a uniformly bounded 
sequence of 
continuous functions, converging pointwise to a function 
which is {\it not\/} 
such a difference. 

In Section 2,  we outline the proof of the main result, 
reducing to a 
qualitative version of the real-variables principle. 
In Section~3, we formulate the Kechris-Louveau invariants, 
which are termed  
here the ``positive transfinite oscillations'', as well as 
some related 
new invariants, called the ``transfinite-oscillations''.
The latter yield a surprising norm-identity on the Banach 
algebra of 
(complex) differences of bounded semi-continuous functions 
(Theorem~3.1). 
Finally, we formulate a quantitative version of the 
real-variables principle  
 (Theorem~3.3), which yields the qualitative version and 
thus the proof of the main result. 
All the results stated here, as well as further 
complements, are given in 
full detail in \cite{R3}. 

\heading%\head 
2.\endheading%. \endhead

Throughout, we use standard Banach space terminology; here 
is a quick review 
of some of the necessary concepts. 
Let $(b_j)$ be a sequence in a Banach space $B$. $(b_j)$ 
is called a 
{\it basic sequence\/} if $(b_j)$ is a basis for its 
closed linear span,  
denoted $[b_j]$. Equivalently, the $b_j$'s are non-zero 
and there is a positive 
$\lambda$ so that $\|\sum_{i=1}^k c_ib_i\| \le 
\lambda\|\sum_{i=1}^n 
c_i b_i\|$ for all $1\le k<n$ and scalars 
$c_1,\ldots,c_n$; if $\lambda$ 
works, we call $(b_j)$ $\lambda$-basic. 
A sequence $(u_j)$ in $B$ is called a {\it block basis\/} 
of $(b_j)$ if 
there exist integers $0\le n_1<n_2<\cdots$ and scalars 
$c_1,c_2,\ldots $ 
so that $u_j= \sum_{i=n_j+1}^{n_{j+1}} c_ib_i$ for all 
$j=1,2,\ldots$; 
$(u_j)$ is called a {\it convex block basis\/} if the 
$c_i$'s can be 
chosen to satisfy $c_i\ge 0$ for all $i$ and $\sum_{i=n_j+
1}^{n_{j+1}}c_i=1$ 
for all $j$. 
Sequences $(x_i)$ and $(y_i)$ in (possibly different) 
Banach spaces are 
called {\it equivalent\/} if there exists an {\it 
isomorphism\/} (i.e., 
a bounded linear invertible) operator $T:[x_i]\to [y_i]$ 
with $Tx_i=y_i$ 
for all $i$. 
If $(b_j)$ is a given basic sequence, $(b_j^*)$ denotes 
the sequence 
of its biorthogonal functionals in $[b_j]^*$, the dual of 
the closed 
linear span of the $b_j$'s. 

We first indicate the reduction of the proof of 
Theorem~1.1 to a ``classical 
real variables'' setting. 
For $K$ a compact metric space, $D(K)$ denotes the set of 
all (complex) 
differences of bounded semi-continuous functions on $K$; 
that is, $f\in D(K)$ if and only if there are bounded lower 
semi-continuous functions $u_1,\ldots,u_4$ on $K$ so that 
$f= (u_1-u_2) 
+ i(u_3-u_4)$: 
equivalently (by results of Baire), if and only if there 
exists a 
sequence $(\varphi_j)$ in $C(K)$ (the space of continuous 
functions on $K$) 
with $\sup_{k\in K} \sum |\varphi_j(k)|<\infty$ and $f=\sum 
\varphi_j$ pointwise. 

The next result follows from refinements of arguments in 
\cite{Bes-P}; 
see Corollary~3.1 of \cite{HOR} (also cf.\ Proposition~1.7 
of \cite{R3}). 

\proclaim{Proposition 2.1} 
Let $K$ be a compact metric space, $f:K\to \complex$ 
discontinuous, 
and $(f_n)$ uniformly bounded in $C(K)$ with $f_n\to f$ 
pointwise. 
Then $f$ is in $D(K)$ if and only if $(f_n)$ has a convex 
block basis 
equivalent to the summing basis. 
\endproclaim 

The $c_0$-theorem then follows immediately from 
Proposition 2.1 and the following 
result (in both Proposition 2.1 and Theorem 2.2, the 
ambient Banach space $B$ is $C(K)$). 

\proclaim{Theorem 2.2} 
Let $K$ be a compact metric space and $(f_n)$ be a 
uniformly bounded 
sequence in $C(K)$ which converges pointwise to a function 
$f$. 
If $f$ is not in $D(K)$, then $(f_n)$ has an $\ss$ 
subsequence. 
\endproclaim

To prove Theorem 2.2, we need some natural permanence 
properties of $\ss$ sequences,  
including the following companion notion. 

\dfn{Definition} 
A basic sequence $(e_j)$ in a Banach space is called {\it 
coefficient 
converging\/}, or $\cc$, if $(\sum_{j=1}^n e_j)$ is a 
weak-Cauchy sequence 
so that whenever scalars $(c_j)$ satisfy $\sup_n 
\|\sum_{j=1}^n c_j e_j\| 
<\infty$, the sequence $(c_j)$ converges. 
\enddfn 

In the sequel, given sequences $(b_j)$ and $(e_j)$ in a 
Banach space, 
$(e_j)$ is called the {\it difference sequence\/} of 
$(b_j)$ if 
$e_1=b_1$ and $e_j=b_j-b_{j-1}$ for all $j>1$. 

\proclaim{Proposition 2.3} 
Let $(b_j)$ be a given sequence in a Banach space. 
%\vskip1pt
\roster
\item"{\rm (a)}" $(b_j)$ is $\ss$ if and only if its 
difference sequence 
is $\cc$.
\item"{\rm (b)}" If $(b_j)$ is $\ss$, every convex block 
basis of $(b_j)$ is 
also $\ss$. 
\item"{\rm (c)}" If $(b_j)$ is a basic sequence, then 
$(b_j)$ is $\ss$ if 
and only if $(b_j^*)$ is $\cc$.
\endroster
\endproclaim 

We note that (c) yields that if $(b_j)$ is $\ss$, 
$[b_j]^*$ is not weakly 
sequentially complete, for it follows easily that if 
$(e_j)$ is $\cc$, 
$(\sum_{j=1}^n e_j)$ is a non-trivial weak Cauchy sequence. 
(b) or (c) combined with 2.1 may be used to show the 
alternatives in the 
$c_0$-theorem are mutually exclusive. 

To obtain $\ss$-sequences, we first note that every 
non-trivial weak-Cauchy 
sequence in a Banach space has an $\s$-subsequence; i.e., 
a weak-Cauchy 
basic sequence which dominates the summing basis (cf.\ 
\cite{HOR}). 
Then we use the following quantitative result, which 
follows by 
diagonalization and an argument of S.~Bellenot \cite{Be}. 

\proclaim{Lemma 2.4} 
Let $(f_j)$ be an $\s$-sequence in a Banach space. 
Then $(f_j)$ has an $\ss$-subsequence provided for every 
$\varep>0$ 
and subsequence $(g_j)$ of $(f_j)$, there is a subsequence 
$(b_j)$ 
of $(g_j)$ with difference sequence $(e_j)$ so that 
whenever $(c_j)$ 
is a sequence of scalars with $c_j=0$ for infinitely many 
$j$ and 
$\|\sum_{j=1}^n c_je_j\|\le 1$ for all $n$,  
$\olim_{j\to\infty} 
|c_j|\le\varep$.
\endproclaim 

We now arrive at the real-variables core of the proof. 

\proclaim{Theorem 2.5}
Let $K$ be a compact metric space and $(f_n)$ be a 
uniformly bounded 
sequence in $C(K)$ which converges pointwise to a function 
$f$ which is 
not in $D(K)$. There exists a $c$ with $|c|=1$ so that 
given $M<\infty$, 
$\kappa >0$, there is a subsequence $(b_j)$ of $(cf_j)$ 
with difference 
sequence $(e_j)$ so that given $0<m_1<m_2<\cdots$ an 
infinite sequence 
of integers, there exists a $t$ in $K$ and an integer $k$ 
with 
%\vskip2pt 
\roster
\item $\sum_{j=1}^k \Re\ e_{m_{2j}} (t) >M$\RM;
%\vskip2pt 
\item
$\Re\ e_{m_{2j}} (t) >0$ for all $1\le j\le k$\RM;
%\vskip2pt 
\item$\sum_{1<i\notin \{m_1,m_2,\dots,m_k\}}  |e_i(t)| 
<\kappa$. 
\endroster
\endproclaim 

\demo{Sketch of proof of Theorem {\rm2.2}, real scalars} 
%\nobreak
We easily reduce to the case where $(f_n)$ is an 
$\s$-sequence satisfying 
the conclusion of Theorem 2.5 with $c=1$. It then follows 
that there are finite 
$\lambda$, $\tau$ so that if $(b_j)$ is a subsequence of 
$(f_j)$ with 
difference sequence $(e_j)$, then 
$$\text{$(e_j)$ is $\lambda$-basic and  $\|e_j^*\| \le 
\tau$\quad for all $j$.}
\leqno(1)$$ 

Now let $\varep >0$, $(f'_j)$ be a subsequence of $(f_j)$, 
$\kappa=1$, 
$M= {1+\tau+\lambda\over\varep}$, and $(b_j)$ be a 
subsequence of $(f'_j)$ 
satisfying the conclusion of Theorem~2.5. 
It follows that the difference sequence $(e_j)$ of $(b_j)$ 
satisfies the conclusion of Lemma~2.4, thus completing the 
proof.\qed 
\enddemo

\heading%\head 
3.\endheading%.\endhead

To obtain Theorem 2.5, we require certain 
transfinite invariants for a general discontinuous 
function.  
We define the {\it transfinite  oscillations\/} of a given 
function 
$f:K\to \complex$; these are similar to invariants 
previously defined in 
\cite{KL}, which we term here the {\it positive 
transfinite oscillations\/}. 
Fix $K$ a compact metric space. 
For $f:K\to [-\infty,\infty]$ an extended real-valued 
function,  
$Uf$ denotes the upper semi-continuous envelope of $f$; 
$Uf(x) = \overline{\lim}_{y\to x} f(y)$ for all $x\in K$. 
(We use non-exclusive lim sups; thus 
$\overline{\lim}_{y\to x} f(y) = 
\inf_U \sup f(U)$, the infimum over all open neighborhoods 
$U$ of $x$.) 

\dfn{Definition}  
Let $f:K\to\complex$ be a given function, $K$ a compact 
metric space, and 
$\alpha$ a countable ordinal. We define the $\alpha^{th}$ 
oscillation of $f$, 
$\osc_\alpha f$, by induction, as follows\,\RM: set 
$\osc_0 f\equiv 0$. 
Suppose $\beta >0$ is a countable ordinal and $\osc_\alpha 
f$ has been 
defined for all $\alpha <\beta$. If $\beta$ is a 
successor, say $\beta = 
\alpha+1$, we define 
$$\widetilde{\osc}_\beta f(x) = \overline{\lim}_{y\to x} 
\bigl( |f(y) - f(x)| 
+ \osc_\alpha f(y)\bigr) \quad\text{for all}\ x\in K\ .$$ 
If $\beta$ is a limit ordinal, we set 
$$\widetilde{\osc}_\beta f = \sup_{\alpha <\beta} 
\osc_\alpha f\ .$$ 
Finally, we set $\osc_\beta f = U \widetilde{\osc}_\beta 
f$. 
\enddfn

The {\it positive $\alpha^{th}$ oscillation\/} of a 
real-valued function $f$, 
$v_\alpha f$, is defined in exactly the same way, except 
that the absolute 
value signs are deleted. 
These are given in \cite{KL}, with a different terminology 
and equivalent 
formulation, where it is established that $f$ 
{\sl is in $D(K)$ if and only 
if $v_\alpha f$ is a bounded function for all $\alpha < 
w_1$. 
Moreover, when this happens, there is an $\alpha <w_1$ 
with $v_\alpha f= 
v_{\alpha+1}f$. Then writing $f= u_\alpha -v_\alpha$, 
$u_\alpha$ is 
upper semi-continuous\/}. 

We obtain a similar result, using the transfinite 
oscillations, which 
yields a rather surprising norm identity. 
We define a norm $\|\cdot\|_D$ on $D(K)$ by 
$$\|f\|_D = \inf \biggl\{ \sup_{k\in K} \sum_j |\varphi_j 
(k)| : 
(\varphi_j) \hbox{ in } C(K)\hbox{ with } \sum \varphi_j = 
f\biggr\}\ .$$ 

\proclaim{Theorem 3.1} 
Let $f:K\to\complex$ be a bounded function. 
Then $f\in D(K)$ if and only if $\osc_\alpha f$ is a 
bounded function for 
all $\alpha <w_1$. 
When this occurs, there is a countable $\alpha$ with 
$\osc_\alpha f= 
\osc_{\alpha+1} f$. If $f$ is real valued, then setting 
$\lambda =\|\,|f| 
+ \osc_\alpha f\|_\infty$, we have that $\|f\|_D=\lambda$, 
and if $u = 
(\lambda-\osc_\alpha f+f)/2$, $v=(\lambda-\osc_\alpha 
f-f)/2$, then $u,v$ 
are non-negative lower semi-continuous functions with 
$f=u-v$ and 
$\|f\|_D = \|u-v\|_\infty$. 
\endproclaim 

It then follows that the ``inf'' in the definition of 
$\|f\|_D$ is actually 
attained (for real-valued $f$). 

Although only the positive transfinite oscillations are 
actually needed for proving the 
$c_0$-theorem, we prefer to use 
the transfinite oscillations, since they  appear more 
natural in  
studying the Banach space structure of $D(K)$ and related 
objects 
(see \cite{R4}). 
The connection between the oscillations is given by the 
following simple 
result. 

\proclaim{Proposition 3.2} 
Let $f:K\to \real$ be a given function and $\alpha $ a 
countable ordinal. Then  
$$v_\alpha (f) \le \osc_\alpha f\le v_\alpha (f) + 
v_\alpha (-f)\ . 
$$ 
\endproclaim 

Theorem 2.5 then follows from Theorem 3.1 and 
Proposition~3.2 (or the cited 
results in \cite{KM}, and the following quantitative 
real-variables result. 

\proclaim{Theorem 3.3} 
Let $(f_j)$ be a uniformly bounded sequence of 
complex-valued  
continuous functions on $K$ a compact metric space, 
converging pointwise 
to a function $f$. Let $\alpha$ be a countable ordinal and 
$x\in K$ be 
given with $0<v_\alpha (\varphi)(x)\dfeq \lambda <\infty$ 
where $\varphi = 
\Re\ f$. Let $\U$ be an open neighborhood of $x$ and 
$0<\eta<1$ be given. 
There exists $(b_j)$ a subsequence of $(f_j)$ with the 
following 
properties\,\RM: Given $1= m_1<m_2<\cdots$ an infinite 
sequence of integers, 
there exist $k$, points $x_1,\ldots, x_{2k-1}$, 
$x_{2k}\dfeq t$ in $\U$,
and positive numbers $\delta_1,\ldots,\delta_k$ so 
that\,\RM: 
\roster%\vskip1pt 
\item
$\varphi (x_{2j}) - \varphi (x_{2j-1}) > (1-\eta)\delta_j$ 
for all $1\le j\le k$,
%\vskip4pt 
\item
$(1+\eta)\lambda > \sum_{j=1}^k \delta_j > 
(1-\eta)\lambda$, 
%\vskip4pt 
\item
$\sum_{m_j\le i<m_{j+1}} |b_i(t) - f(x_j)| <\eta 
\delta_{[{j+1\over2}]} $ for all $1\le j\le 2k-1$,
%\vskipq4pt 
\item
$\sum_{i\ge m_{2k}} |b_i(t) - f(t)| <\eta\, \delta_k$\.
\endroster
\endproclaim 

The proof of Theorem 3.3 is given by transfinite 
induction, and the 
formulation in terms of arbitrary open neighborhoods $\U$ 
is crucial, 
although this is not used in the proof of Theorem~2.5. 
The argument is easily reduced to the ``$\alpha$ to 
$\alpha +1$'' case. 

\demo{Sketch of proof of the induction step} 
We suppose the result proved for $\alpha$ and assume 
$0<v_\alpha(\varphi)(x) 
<v_{\alpha+1} (\varphi)(x) \dfeq \beta$. 
Now given $\eta>0$ and $\U$ an open neighborhood of $x$, 
we obtain the 
existence of positive numbers $\bar\lambda$ and $\delta$ 
with 
$(1-\eta)\beta <\bar\lambda +\delta <(1+\eta)\beta$, 
$x_1\in\U$, and a 
subsequence $(b_j)$ of $(f_j)$, so that given $1<r\le s$, 
there is an 
open set $\V\subset \U$ and an $x_2\in\V$ with 
$$\gather
\varphi (x_2)-\varphi(x_1) > (1-\eta)\delta; \tag 2\\ 
\sum_{1\le i<r}|b_i(t)- f(x_1)| <\eta\delta\ \text{ for 
all }\ t\in\V;\tag 3\\
\sum_{r\le i\le s} |b_i(t) - f(x_2)| <\eta\delta\ \text{ 
for all }\ t\in \V;
\tag 4\\
\un{\lambda} \le \lambda < (1+\eta) \beta-\delta\ ,\ 
\text{ where }\ 
\lambda = v_\alpha (\varphi)(x_2)\. \tag 5 
\endgather$$ 
The numbers $\bar\lambda$ and $\delta$ are obtained using 
the definition 
of $v_{\alpha+1}(\varphi)$, and $(b_i)$ is obtained from 
the proof of the 
``$\alpha=1$'' case. 
Finally, after a further refinement, using the inductive 
hypothesis, we obtain 
$$\left\{ \eqalign{&\text{$(b_i)_{i>s}$ satisfies the 
conclusion of Theorem 
3.3}\cr 
&\text{for the $\alpha$-case, with ``$\U$'' $=\V$, ``$x$'' 
$=x_2$.}\cr} 
\right. \leqno(6)$$ 

Now thanks to (2), (5) and the fact that the inequalities 
in (3), (4) hold for 
{\it all\/} $t\in \V$, we obtain that $(b_i)$ satisfies the 
conclusion of Theorem~3.3 for the $(\alpha+1)$-case.\qed
\enddemo

\enddocument